\documentclass[11pt]{article}
\usepackage{exscale,relsize}
\usepackage{amsmath}
\usepackage{amsfonts}
\usepackage{amssymb}
\usepackage{calc}
\usepackage{theorem}
\usepackage{pifont}      %needed by dingautolist
\newcommand{\GX}{\ensuremath{\Gamma}}

\newcommand{\scal}[2]{\langle{{#1},{#2}}\rangle}

\newcommand{\yosida}{\ensuremath{ \; {}^}}

\newcommand{\RR}{\ensuremath{\mathbb R}}
\newcommand{\SN}{\ensuremath{\mathbb S}^N}
\newcommand{\SNP}{\ensuremath{\mathbb S}^N_+}
\newcommand{\SNPP}{\ensuremath{\mathbb S}^N_{++}}

\newcommand{\RPP}{\ensuremath{\,\left]0,+\infty\right[}}
\newcommand{\RX}{\ensuremath{\,\left]-\infty,+\infty\right]}}

\newcommand{\thalb}{\ensuremath{\tfrac{1}{2}}}

\newcommand{\fettx}{\ensuremath{\boldsymbol{x}}}
\newcommand{\fettla}{\ensuremath{\boldsymbol{\lambda}}}
\newcommand{\bA}{\ensuremath{\boldsymbol{A}}}
\newcommand{\average}{\ensuremath{\mathcal{R}_{\mu}(\bA,\fettla)}}

\newcommand{\averageonelambda}{\ensuremath{\mathcal{R}_{1}(\bA,\fettla)}}
\newcommand{\averagex}{\ensuremath{\mathcal{R}_{}(\fettx,\fettla)}}
\newcommand{\fetty}{\ensuremath{\boldsymbol{y}}}
\newcommand{\averagey}{\ensuremath{\mathcal{R}_{}(\fetty,\fettla)}}
\newcommand{\averageinversex}{\ensuremath{\mathcal{R}_{}(\fettx^{-1},\fettla)}}
\newcommand{\averageoneinverse}{\ensuremath{\mathcal{R}_{1}(\bA^{-1},\fettla)}}
\newcommand{\geoaverage}{\ensuremath{\mathcal{G}_{}(\fettx,\fettla)}}
\newcommand{\geoaveragey}{\ensuremath{\mathcal{G}_{}(\fetty,\fettla)}}
\newcommand{\geoinverse}{\ensuremath{\mathcal{G}_{}(\fettx^{-1},\fettla)}}
\newcommand{\averageinverse}{\ensuremath{\mathcal{R}_{\mu^{-1}}(\bA^{-1},\fettla)}}

\newcommand{\res}{\ensuremath{\mathcal{R}_{\mu}}}

\newcommand{\harm}{\ensuremath{\mathcal{H}(\bA,\fettla)}}
\newcommand{\arithmetic}{\ensuremath{\mathcal{A}(\bA,\fettla)}}
\newcommand{\harminverse}{\ensuremath{\mathcal{H}(\bA^{-1},\fettla)}}
\newcommand{\arithmeticinverse}{\ensuremath{\mathcal{A}(\bA^{-1},\fettla)}}
\newcommand{\bB}{\ensuremath{\boldsymbol{B}}}
\newcommand{\averageb}{\ensuremath{\mathcal{R}_{\mu}(\bB,\fettla)}}

\newcommand{\pflm}{\ensuremath{p_\mu({\fettf,\fettla)}}}

\newcommand{\Id}{\ensuremath{\operatorname{Id}}}

\newcommand{\fettf}{\ensuremath{\boldsymbol{f}}}

\newcommand{\pinf}{\ensuremath{+\infty}}

\renewcommand{\phi}{\ensuremath{\varphi}}
\newtheorem{theorem}{Theorem}[section]
\newtheorem{lemma}[theorem]{Lemma}
\newtheorem{fact}[theorem]{Fact}
\newtheorem{corollary}[theorem]{Corollary}
\newtheorem{proposition}[theorem]{Proposition}
\newtheorem{definition}[theorem]{Definition}

\theoremstyle{plain}{\theorembodyfont{\rmfamily}
}
\theoremstyle{plain}{\theorembodyfont{\rmfamily}
}
\theoremstyle{plain}{\theorembodyfont{\rmfamily}
}
\theoremstyle{plain}{\theorembodyfont{\rmfamily}
\newtheorem{example}[theorem]{Example}}
\theoremstyle{plain}{\theorembodyfont{\rmfamily}
\newtheorem{remark}[theorem]{Remark}}
\theoremstyle{plain}{\theorembodyfont{\rmfamily}
}

\newcommand{\pluss}{{\hskip1pt \raise1pt\vbox{\hrule width6pt \vskip1pt
\hrule width6pt}\kern-4pt{\lower1pt\hbox{\vrule height6pt \kern1pt\vrule
height6pt}}\hskip5pt}}
\newcommand{\timess}{\star}

\begin{document}

\title{{\sffamily The Resolvent Average for Positive Semidefinite Matrices}}

\author{
Heinz H.\ Bauschke\thanks{Mathematics, Irving K.\ Barber School,
University of British Columbia Okanagan, Kelowna, British Columbia
V1V 1V7, Canada. E-mail: \texttt{heinz.bauschke@ubc.ca}.},~
Sarah M.\ Moffat\thanks{Mathematics, Irving K.\ Barber School, University
of British Columbia Okanagan, Kelowna, British Columbia V1V 1V7,
Canada. E-mail: \texttt{smoffat99@gmail.com}.}~ and
Xianfu Wang\thanks{Mathematics, Irving K.\ Barber School, University
of British Columbia Okanagan, Kelowna, British Columbia V1V 1V7,
Canada. E-mail: \texttt{shawn.wang@ubc.ca}.}}

\date{October 19, 2009}

\maketitle

% \vskip 8mm

\begin{abstract} \noindent
We define a new average --- termed the \emph{resolvent average} ---
for positive semidefinite matrices.
For positive definite matrices, the resolvent average
enjoys self-duality and it interpolates between
the harmonic and the arithmetic averages, which it approaches when
taking appropriate limits.
We compare the resolvent average to the geometric mean.
Some applications to matrix functions are also given.
\end{abstract}



%%%%%%%%%%%%%%%%%%%%%%%%%%%%%%%%%%%%%%%%%%%%%%%%%%%%%%%%%%%%%%%%%%%%%%%%%%%%%%%%%%%%%%%%%%%%%%%%
%%%%%%%%%%%%%%%%%%%%%%%%%%%%%%%%%%%%%%%%%%%%%%%%%%%%%%%%%%%%%%%%%%%%%%%%%%%%%%%%%%%%%%%%%%%%%%%%
%%%%%%%%%%%%%%%%%%%%%%%%%%%%%%%%%%%%%%%%%%%%%%%%%%%%%%%%%%%%%%%%%%%%%%%%%%%%%%%%%%%%%%%%%%%%%%%%
\section{Introduction}

Let $A_{i}, i=1,\ldots,n$ be positive semidefinite matrices, $\lambda_{i}>0$
with $\sum_{i=1}^{n}\lambda_{i}=1$ and $\Id:\RR^N\rightarrow\RR^N$ be the identity mapping. For
$$\bA=(A_{1},\ldots,A_{n}),\quad
\fettla=(\lambda_{1},\ldots, \lambda_{n}),$$ we define
\begin{equation}\label{thestart}
\average =\big[\lambda_{1}(A_{1}+\mu^{-1}\Id)^{-1}+\cdots+\lambda_{n}(A_{n}+\mu^{-1}\Id)^{-1}\big]^{-1}
-\mu^{-1}\Id,
\end{equation}
and call it the \emph{resolvent average} of $\bA$. This is motivated
from the fact that when $\mu=1$
\begin{equation}\label{original}
\big(\averageonelambda+\Id\big)^{-1}=\lambda_{1}\big(A_{1}+\Id\big)^{-1}+\cdots+\lambda_{n}\big(
A_{n}+\Id\big)^{-1},
\end{equation}
which says that the resolvent of $\averageonelambda$ is the (arithmetic)
average of resolvents of the $A_{i}$, with weight
$\fettla=(\lambda_{1},\ldots,\lambda_{n})$. The resolvent average provides
a novel averaging technique, and having the parameter $\mu$
in $\average$ will allow us to take limits later on.
We denote the well known
\emph{harmonic average} and \emph{arithmetic average} by
$$\harm=(\lambda_{1}A_{1}^{-1}+\cdots+\lambda_{n}A_{n}^{-1})^{-1},$$
$$\arithmetic=\lambda_{1}A_{1}+\cdots+\lambda_{n}A_{n},$$
respectively.
In the literature, $(A_{1}^{-1}+\cdots +A_{n}^{-1})^{-1}$ is called the
\emph{parallel sum} of the matrices $A_{1},\ldots,
A_{n}$; see, e.g.,
\cite{anderson,attouch,urruty1,urruty2,mazure,Mou1,Mou2,passty}.

\emph{The goal of this note is to study relationships among the
resolvent average, the harmonic average and the arithmetic average
of matrices. Our proofs are based on convex analytical techniques
and on the proximal average, instead of the more commonly employed
matrix diagonalizations.}

The plan of the paper is as follows. After proving some elementary
properties of $\average$ in Section~\ref{basic}, we gather some basic
properties of proximal averages and general convex functions in
Section~\ref{keyfacts}. The main results, which are given
in Section~\ref{takinglimit}, state that
$$\harm\preceq\average\preceq\arithmetic,$$
$$\lim_{\mu\rightarrow 0^+}\average=\arithmetic, \quad \lim_{\mu\rightarrow +\infty}\average=\harm,$$
and that $\average$ enjoys self-duality, namely $\big[\average\big]^{-1}=\averageinverse.$
In Section~\ref{compare}, we show that the resolvent average and geometric
mean have strikingly similar properties, even though they are different.

\noindent {\bf Notation}:
Throughout, $\RR^N$ is the standard $N$-dimensional Euclidean space.
For $\lambda >0$,
\begin{equation}\label{regularization}
J_{A}=(\Id+A)^{-1}, \quad \yosida{\lambda}A=\lambda^{-1}(\Id-J_{\lambda A}),
\end{equation}
are called the \emph{resolvent} of $A$ and \emph{Yosida $\lambda$-regularization } of $A$. A function
$f:\RR^N\rightarrow\RX = \RR\cup\{\pinf\}$
is said to be convex if its domain is convex and
\begin{equation}\label{convex}
f(\lambda x+(1-\lambda)y)\leq \lambda f(x)+(1-\lambda)f(y)\quad \forall \ x,y \in\RR^N, 0<\lambda <1,
\end{equation}
with $f$ being strictly convex if \eqref{convex} becomes a strict inequality whenever $x\neq y$.
The function $f$ is proper
if $f(x)>-\infty\ \forall x\in\RR^N$ and $f(x_{0})<\pinf$
for some $x_{0}\in\RR^N$.
The class of proper
lower semicontinuous convex functions from $\RR^N\rightarrow\RX$ will be denoted by $\GX$. For $f\in\GX$, $\partial f$ denotes its convex subdifferential:
$\partial f(x)=\{x^*\in\RR^N:\ f(y)\geq f(x)+\scal{x^*}{y-x}\ \forall y\in \RR^N\}.$
If $f$ is differentiable at $x$, then $\partial f(x)=\{\nabla f(x)\}$.
$f^*$ denotes its \emph{Fenchel conjugate} given by $(\forall x^*\in\RR^N)\ f^*(x^*)=\sup_{x}\{\scal{x^*}{x}-f(x)\}$.
For $\alpha>0$, $\alpha\timess f=\alpha f(\cdot/\alpha)$. If $f,g\in\GX$, $f\Box g$ stands for the infimal convolution
of $f,g$ given by
$(f\Box g)(x)
=\inf\{f(x_{1})+g(x_{2}):\ x_{1}+x_2=x\}$ $\forall x\in \RR^N$.
When $A:\RR^N\rightarrow\RR^N$ is linear, the quadratic form $q_{A}:\RR^N\rightarrow\RR$ is defined by
$$q_{A}(x)=\thalb\scal{Ax}{x} \ \forall x\in\RR^N,$$
 and we also use $q_{\Id}=j$ interchangeably. For convex functions $f_{1},\ldots, f_{n}$, we write
 $$\fettf=(f_{1},\ldots, f_{n}), \qquad \fettf^*=(f_{1}^*,\ldots, f_{n}^*).$$
In the space $\SN$ of $N\times N$ real symmetric matrices,
$\SNP$ (resp. $\SNPP$) denotes the set of $N\times N$ positive semidefinite matrices (resp. positive definite matrices).
 For $X,Y\in \SN$, we write $Y\preceq X$  if $X-Y\in \SNP$ and
 $Y\prec X$ if $X-Y\in \SNPP$.

\section{Basic properties}\label{basic}

In this section, we give some basic properties of $\average$.
\begin{proposition}\label{differentreform}
We have
\begin{equation}\label{resolventidentity}
J_{\mu \average}=\lambda_{1}J_{\mu A_{1}}+\cdots +\lambda_{n}J_{\mu A_{n}},
\end{equation}
\begin{equation}\label{t:averageyosida}
\yosida{\mu}\big(\average\big)=\lambda_{1}\yosida{\mu}A_{1}+\cdots +\lambda_{n}\yosida{\mu}A_{n}.
\end{equation}
\end{proposition}
\begin{proof}
Multiplying \eqref{thestart} both sides by $\mu$ gives
$$
\mu\average +\Id =\big[\lambda_{1}(\mu A_{1}+\Id)^{-1}+\cdots+\lambda_{n}(\mu A_{n}+\Id)^{-1}\big]^{-1}.
$$
Then \eqref{resolventidentity} follows by taking inverse both sides and using \eqref{regularization}.

By \eqref{resolventidentity}, we obtain that
$$(\Id-J_{\mu \average})=\lambda_{1}(\Id -J_{\mu A_{1}})+\cdots +\lambda_{n}(\Id -J_{\mu A_{n}}).$$
Dividing both sides by $\mu$,
$$\mu^{-1}(\Id-J_{\mu \average})=\lambda_{1}\mu^{-1}(\Id -J_{\mu A_{1}})+\cdots +\lambda_{n}\mu^{-1}(\Id -J_{\mu A_{n}}).$$
It remains to use \eqref{regularization}.
\end{proof}

\begin{proposition}\label{primalinverse}
Let $\bA=(A_{1},A_{1}^{-1},\ldots, A_{m}, A_{m}^{-1})$,
$\lambda=(\tfrac{1}{2m},\tfrac{1}{2m},\ldots, \tfrac{1}{2m})$,
and $\mu=1$.
Then $\average =\Id$.
\end{proposition}
\begin{proof} This follows from \eqref{original} and the
identity
$(A+\Id)^{-1}+(A^{-1}+\Id)^{-1}=\Id.$
\end{proof}
\begin{proposition}\label{sameoperator} Let $\bA=(A_{1},\ldots, A_{1})$.
Then $\average =A_{1}$.
\end{proposition}
\begin{proof}
 We have
\begin{align*}
\average & =\big((\lambda_{1}+\cdots+\lambda_{n})(A_{1}+\mu^{-1}\Id)^{-1}\big)^{-1}-\mu^{-1}\Id\\
& =
\big((A_{1}+\mu^{-1}\Id)^{-1}\big)^{-1}-\mu\Id=
A_{1}+\mu^{-1}\Id-\mu^{-1}\Id=A_{1},
\end{align*}
which proves the result.
\end{proof}

Note that for $A,B\in \SNPP$, we have
\begin{equation}\label{comparenew}
A\succeq B\quad \Leftrightarrow \quad A^{-1}\preceq B^{-1}
\end{equation}
and
\begin{equation}\label{compareneu}
A\succ B\quad \Leftrightarrow \quad A^{-1}\prec B^{-1};
\end{equation}
see, e.g., \cite[Corollary~7.7.4.(a)]{HoJo} and
\cite[Section~16.E]{MO} or \cite[page~55]{borwein}.

\begin{proposition}\label{monotonecompare}
 Assume that $(\forall\ i)\ A_{i}, B_{i}\in \SNP$ and $A_{i}\succeq B_{i}$.
Then
\begin{equation}\label{monotoneorder}
\average\succeq \averageb.
\end{equation}
Furthermore, if additionally some $A_j\succ B_j$, then $\average\succ \averageb$.
\end{proposition}
\begin{proof}
Note that $\forall \ \mu>0$,
$$A_{i}+\mu^{-1}\Id\succeq B_{i}+\mu^{-1}\Id\succ 0,$$
so that
$$0 \prec (A_{i}+\mu^{-1}\Id)^{-1}\preceq (B_{i}+\mu^{-1}\Id)^{-1},$$
by \eqref{comparenew}. As $\SNP$ and $\SNPP$
are convex cones, we obtain that
\begin{equation}\label{keeppositive}
0\prec \sum_{i=1}^{n}\lambda_{i}(A_{i}+\mu^{-1}\Id)^{-1}\preceq \sum_{i=1}^{n}\lambda_{i}(B_{i}+\mu^{-1}\Id)^{-1}.
\end{equation}
Using \eqref{comparenew} on \eqref{keeppositive}, followed by subtracting $\mu^{-1}\Id$,  gives
$$\big[\sum_{i=1}^{n}\lambda_{i}(A_{i}+\mu^{-1}\Id)^{-1}\big]^{-1}-\mu^{-1}\Id\succeq \big[\sum_{i=1}^{n}\lambda_{i}(B_{i}+\mu^{-1}\Id)^{-1}\big]^{-1}-\mu^{-1}\Id,$$
which establishes \eqref{monotoneorder}.
The ``Furthermore'' part follows analogously using \eqref{compareneu}.
\end{proof}

\begin{theorem}
\label{t:blabla}
Assume that $(\forall\ i)\ A_{i}\in \SNP$.
Then $\average\in \SNP$.
Furthermore, if additionally some $A_j\in\SNPP$, then
$\average\in \SNPP$.
\end{theorem}
\begin{proof}
This follows from Proposition~\ref{monotonecompare} (with each $B_i=0$) and
Proposition~\ref{sameoperator}.
\end{proof}

We end this section with a recursion formula that may be
verified directly using the definitions.

\begin{proposition}[recursion]
We have
$$\res(A_{1},\ldots,A_n;\lambda_{1},\ldots, \lambda_{n})=
\res\Big(\res\big(A_1,\ldots,A_{n-1};\tfrac{\lambda_1}{1-\lambda_n},
\ldots,\tfrac{\lambda_{n-1}}{1-\lambda_n}\big),
A_n;1-\lambda_n,\lambda_n\Big).$$
\end{proposition}

\section{Auxiliary results and facts}\label{keyfacts}
The key tool in this note is the \emph{proximal average of convex
functions}, which finds its roots in \cite{matou,Minty,Moreau},
and which has been further systematically studied in
\cite{bwang,bglw,BLT,byw}.

\begin{definition}[proximal average] \label{d:proxaverage} Let $(\forall i) \ f_{i}\in\GX$.
The $\lambda$-weighted proximal average of $\fettf=(f_1, \ldots, f_{n})$ with parameter $\mu$ is
defined by
\begin{equation}
\pflm=
\bigg(\lambda_{1} (f_1+\tfrac{1}{\mu}j)^* +
\lambda_{2}(f_2+\tfrac{1}{\mu}j)^*+\cdots +\lambda_{n}(f_{n}+\tfrac{1}{\mu}j)^{*}\bigg)^* -\tfrac{1}{\mu}
j.
\end{equation}
\end{definition}
The function
$\pflm$ is a proper lower semicontinuous convex function on $\RR^N$, and
it inherits many desirable properties from each underlying function
$f_{i}$; see \cite{bglw,BLT}.
A fundamental property of proximal average is:
\begin{fact}\emph{(\cite[Theorem 5.1]{bglw})}\label{t:Fenchel}
$\big(\pflm\big)^* = p_{\mu^{-1}}({\fettf^*,\fettla})$.
\end{fact}

To give new proofs of Fact~\ref{t:ineqs} and Fact~\ref{t:mono} below,
we shall need reformulations of $\pflm$.

\begin{proposition}\label{domapp1} Let $f_{1},\ldots,
f_{n}\in \GX$ and $\lambda_{1},\ldots, \lambda_{n}> 0$ with
$\sum_{i=1}^{n}\lambda_{i}=1$.
Then for every $x\in \RR^n$, \\
$\pflm(x)$
\begin{align}
&=\inf_{x_{1}+\cdots+x_{n}=x}\left\{\lambda_{1}(f_{1}+\frac{1}{\mu}j)(\frac{x_{1}}{\lambda_{1}})+\cdots
+\lambda_{n}(f_{n}+\frac{1}{\mu}j)(\frac{x_{n}}{\lambda_{n}})\right\}-\frac{1}{\mu}j(x)\label{onerep}\\
&=\inf_{x_{1}+\cdots+x_{n}=x}\left\{\lambda_{1}f_{1}(\frac{x_{1}}{\lambda_{1}})+\cdots
+\lambda_{n}f_{n}(\frac{x_{n}}{\lambda_{n}})+\frac{1}{4\mu}\sum_{i=1}^{n}
\sum_{j=1}^{n}\lambda_{i}\lambda_{j}\|\frac{x_{i}}{\lambda_{i}}-
\frac{x_{j}}{\lambda_{j}}\|^2\right\}\label{tworep}\\
&=\inf_{\lambda_{1}y_1+\cdots +\lambda_{n}y_{n}=x}
\left\{\lambda_{1} f_{1}(y_{1})+\cdots
+\lambda_{n}f_{n}(y_{n})+\frac{1}{\mu}[\lambda_{1}j(y_{1})+\cdots
+\lambda_{n}j(y_{n})-j(\lambda_{1}y_{1}+\cdots
+\lambda_{n}y_{n})]\right\}\label{fourrep}\\
&= \inf_{\lambda_{1}y_1+\cdots
+\lambda_{n}y_{n}=x}\left\{\lambda_{1}f_{1}(y_{1})+\cdots+\lambda_{n}f_{n}(y_{n})+
\frac{1}{4\mu}\sum_{i=1}^{n}\sum_{j=1}^{n}
\lambda_{i}\lambda_{j}\|y_{i}-y_{j}\|^2\right\}\label{threerep}\\
&=\inf_{x_{1}+\cdots
+x_{n}=x}\left\{\lambda_{1}f_{1}(\frac{x_{1}}{\lambda_{1}})+\cdots
+\lambda_{n}f_{n}(\frac{x_{n}}{\lambda_{n}})+
\frac{1}{\mu}[\lambda_{1}j(x-\frac{x_{1}}{\lambda_{1}}))+
\cdots +\lambda_{n}j(x-\frac{x_{n}}{\lambda_{n}})]\right\}.\label{fiverep}
\end{align}
Furthermore, the infimal convolutions in
(\ref{onerep})--(\ref{fiverep}) are exact.
\end{proposition}
\begin{proof} Indeed, as
$$\left(f_{i}+\frac{1}{\mu}j\right)^*=f_{i}^*\Box (\mu j),$$
it is finite-valued everywhere, we write
$$
f=\lambda_{1}\star(f_{1}+\frac{1}{\mu}j)\Box \cdots \Box
\lambda_{n}\star(f_{n}+\frac{1}{\mu}j)- \frac{1}{\mu}j,$$ by
\cite[Theorem~16.4]{Rocky}. That is, for every $x$,
$$f(x)=\inf\left\{\lambda_{1}(f_{1}+\frac{1}{\mu}j)(\frac{x_{1}}{\lambda_{1}})
+\cdots
+\lambda_{n}(f_{n}+\frac{1}{\mu}j)(\frac{x_{n}}{\lambda_{n}}): \;
x_{1}+\cdots +x_{n}=x\right\}-\frac{1}{\mu}j(x),$$ and the infimum
is attained. Hence (\ref{onerep}) holds.

Now rewrite (\ref{onerep}) as
\begin{align}
&\inf_{x_{1}+\cdots+x_{n}=x}\left\{\lambda_{1}f_{1}(\frac{x_{1}}{\lambda_{1}})+\cdots
+\lambda_{n}f_{n}(\frac{x_{n}}{\lambda_{n}})+\frac{1}{\mu}[\lambda_{1}j(\frac{x_{1}}{\lambda_{1}})
+\cdots +\lambda_{n}j(\frac{x_{n}}{\lambda_{n}})-j(x_{1}+\cdots+
x_{n})]\right\},\label{whom}\\
&=\inf_{\lambda_{1}y_1+\cdots +\lambda_{n}y_{n}=x}
\left\{\lambda_{1} f_{1}(y_{1})+\cdots
+\lambda_{n}f_{n}(y_{n})+\frac{1}{\mu}[\lambda_{1}j(y_{1})+\cdots
+\lambda_{n}j(y_{n})-j(\lambda_{1}y_{1}+\cdots
+\lambda_{n}y_{n})]\right\}.\nonumber
\end{align}
Thus, (\ref{tworep})--(\ref{threerep}) follow by using the
identity
$$\sum_{i=1}^{n}\lambda_{i}j(y_{i})-j(\sum_{i=1}^{n}\lambda_{i}y_{i})=
\frac{1}{4}\sum_{i=1}^{n}\sum_{j=1}^{n}\lambda_{i}\lambda_{j}\|y_{i}-y_{j}\|^2.$$
Observe that
\begin{align*}
&\lambda_{1}j(x_{1}+\cdots
+x_{n}-\frac{x_{1}}{\lambda_{1}})+\cdots
+\lambda_{n}j(x_{1}+\cdots+x_{n}-\frac{x_{n}}{\lambda_{n}})\\
&=\lambda_{1}j(\frac{x_{1}}{\lambda_{1}})+\cdots+\lambda_{n}j(\frac{x_{n}}{\lambda_{n}})-j(x_{1}+\cdots+x_{n}),
\end{align*}
we have (\ref{fiverep}) by (\ref{whom}).
\end{proof}

\begin{fact}\emph{(\cite[Theorem 5.4]{bglw})}\label{t:ineqs}
$(\lambda_1f_1^* + \cdots + \lambda_n f_n^*)^{*} \leq \pflm
\leq \lambda_1 f_1+\cdots +\lambda_n f_n$.
\end{fact}
\begin{proof}
This follows from \eqref{tworep} or \eqref{threerep}.
\end{proof}

\begin{fact}\emph{(\cite[Theorem 8.5]{bglw})}\label{t:mono}
Let $x\in \RR^N$. Then the function
\begin{equation} \label{e:hike:0}
\RPP\to\RX\colon\mu \mapsto \pflm(x)\quad\text{is decreasing.}
\end{equation}
Consequently, $\lim_{\mu\to 0^+}\pflm(x)$ and $\lim_{\mu\to\pinf}\pflm(x)$ exist. In fact,
\begin{equation} \label{e:hike:a}
\lim_{\mu\to 0^+}\pflm(x) = \sup_{\mu>0}\pflm(x)=
\big(\lambda_1f_1+\cdots+\lambda_nf_n\big)(x)
\end{equation}
and
\begin{equation} \label{e:hike:b}
\lim_{\mu\to \pinf}\pflm(x) =\inf_{\mu>0} \pflm(x)
= \big(\lambda_1\timess f_1\Box
\cdots\Box\lambda_n\timess f_n\big)(x).
\end{equation}
\end{fact}
\begin{proof}
\eqref{e:hike:0} follows from \eqref{tworep}. \eqref{e:hike:b} also follows from \eqref{tworep}.

To see \eqref{e:hike:a},
by \eqref{fiverep}, $\forall x\in \RR^N$,
\begin{align*}
\pflm(x) &\geq
\lambda_{1}\inf_{x_{1}}\bigg(f_{1}(x_{1}/\lambda_{1})+\frac{1}{\mu}j(x-x_{1}/\lambda_{1})\bigg)
+\cdots
+\lambda_{n}\inf_{x_{n}}\bigg(f_{n}(x_{n}/\lambda_{n})+\frac{1}{\mu}j(x-x_{n}/\lambda_{n})\bigg)\\
&=\lambda_{1}e_{\mu}f_{1}(x)+\cdots+\lambda_{n}e_{\mu}f_{n}(x),
\end{align*}
where $e_{\mu}f_{i}=f_{i}\Box (1/\mu j)$.
Then
$$ \lambda_{1}e_{\mu}f_{1}+\cdots+\lambda_{n}e_{\mu}f_{n}\leq \pflm \leq
\lambda_{1}f_{1}+\cdots+\lambda_{n}f_{n},
$$
so that
$$\lim_{\mu\rightarrow 0^{+}}\pflm=\lambda_{1}f_{1}+\cdots+\lambda_{n}f_{n},$$
since $\lim_{\mu\rightarrow 0^+} e_{\mu}f_{i}=f_{i}$
by \cite[Theorem~2.26 and Theorem~1.25]{RockWets}.
\end{proof}

\begin{fact}\emph{(\cite[Theorem 25.7]{Rocky})}\label{rockstar}
Let $C$ be a nonempty open convex subset of $\RR^N$,
and let $f$ be a convex function
which is finite and differentiable on $C$. Let $f_{1},f_{2},\ldots,$ be a sequence of convex functions finite
and differentiable on $C$ such that $\lim_{i\rightarrow\infty} f_{i}(x)=f(x)$ for every
$x\in C$. Then $$\lim_{i}\nabla f_{i}(x)=\nabla f(x), \quad \forall x\in C.$$
In fact, the sequence of gradients
$\nabla f_{i}$ converges to $\nabla f$ uniformly on every compact
subset of $C$.
\end{fact}

\begin{fact}\emph{(\cite[page 108]{Rocky})}\label{qconjugate}
Let $Q\in \SNPP$. Then $(q_{Q})^*=q_{Q^{-1}}.$
\end{fact}

\begin{fact}\emph{(\cite[Theorem 23.5]{Rocky})}\label{conjgrad}
Let $f:\RR^N\rightarrow\RX$ be a proper lower semicontinuous convex function. Then $\partial f^*=(\partial f)^{-1}.$
\end{fact}

\section{Main results}\label{takinglimit}

We start by computing the proximal average of general linear-quadratic
functions thereby extending
\cite[Example~4.5]{bglw} and \cite[Example~7.4]{BLT}.

\begin{lemma}\label{quadraticprox}
Let $A_{i}\in \SNP$, $b_{i}\in\RR^N$, $r_{i}\in \RR$. If each $f_{i}=q_{A_{i}}+\scal{b_{i}}{\cdot}+r_{i}$, i.e.,
linear-quadratic,
then $\forall x^*$,
\begin{align}\label{quadratic}
\pflm(x^*)& =q_{\average}(x^*)+\scal{x^*}{\big(\sum_{i=1}^{n}\lambda_{i}(A_{i}+\mu^{-1}\Id)^{-1}\big)^{-1}\sum_{i=1}^{n}
\lambda_{i}(A_{i}+\mu^{-1}\Id)^{-1}b_{i}}
+\nonumber\\
& q_{\big(\sum_{i=1}^{n}\lambda_{i}(A_{i}+\mu^{-1}\Id)^{-1}\big)^{-1}}(\sum_{i=1}^{n}
\lambda_{i}(A_{i}+\mu^{-1}\Id)^{-1}b_{i})-\sum_{i=1}^{n}\lambda_{i}\big(q_{(A_{i}+\mu^{-1}\Id)^{-1}}(b_{i})-r_{i}\big).
\end{align}
In particular,
if $(\forall i)\ f_{i}$ is quadratic, i.e., $\ b_{i}=0, r_i=0$, then $\pflm$ is quadratic with
$$\pflm =q_{\average};$$
If $(\forall i)\ f_{i}$ is affine, i.e., $\ A_{i}=0$, then
$\pflm$ is affine.
\end{lemma}
\begin{proof}
We have
$f_{i}+\mu^{-1}j=q_{(A_{i}+\mu^{-1}\Id)}+\scal{b_{i}}{\cdot}+r_{i}$ and by Fact~\ref{qconjugate}
\begin{align*}
(f_{i}+\mu^{-1}j)^*(x^*)& =q_{(A_{i}+\mu^{-1}\Id)^{-1}}(x^*-b_{i})-r_{i}\\
&=q_{(A_{i}+\mu^{-1}\Id)^{-1}}(x^*)-\scal{x^*}{(A_{i}+\mu^{-1}\Id)^{-1}b_{i}}+q_{(A_{i}+\mu^{-1}\Id)^{-1}}(b_{i})-
r_{i}.
\end{align*}
Then
$\big(\lambda_{1}(f_{1}+\mu^{-1}j)^{*}+\cdots+\lambda_{n}(f_{n}+\mu^{-1}j)^*\big)(x^*)=$
\begin{align*}
&\sum_{i=1}^{n}\lambda_{i}\bigg(q_{(A_{i}+\mu^{-1}\Id)^{-1}}(x^*)-\scal{x^*}{(A_{i}+\mu^{-1}\Id)^{-1}b_{i}}+q_{(A_{i}+\mu^{-1}\Id)^{-1}}(b_{i})-
r_{i}\bigg)
\\
&= q_{{\sum_{i=1}^{n}\lambda_{i}(A_{i}+\mu^{-1}\Id)^{-1}}}(x^*)-
\scal{x^*}{\sum_{i=1}^{n}\lambda_{i}(A_{i}+\mu^{-1}\Id)^{-1}b_{i}}
+\sum_{i=1}^{n}\lambda_{i}\big(q_{(A_{i}+\mu^{-1}\Id)^{-1}}(b_{i})-
r_{i}\big).
\end{align*}
It follows that
$\pflm(x^*)=$
$$q_{[{\sum_{i=1}^{n}\lambda_{i}(A_{i}+\mu^{-1}\Id)^{-1}}]^{-1}}(x^*+\sum_{i=1}^{n}\lambda_{i}(A_{i}+\mu^{-1}\Id)^{-1}b_{i})
-\sum_{i=1}^{n}\lambda_{i}\big(q_{(A_{i}+\mu^{-1}\Id)^{-1}}(b_{i})-
r_{i}\big)-q_{\mu^{-1}\Id}(x^*).$$
As
\begin{align*}
& q_{{[{\sum_{i=1}^{n}\lambda_{i}(A_{i}+\mu^{-1}\Id)^{-1}}]^{-1}}}(x^*+
\sum_{i=1}^{n}\lambda_{i}(A_{i}+\mu^{-1}\Id)^{-1}b_{i})
\\
&=q_{[{\sum_{i=1}^{n}\lambda_{i}(A_{i}+\mu^{-1}\Id)^{-1}}]^{-1}}(x^*)
+\scal{x^*}{[\sum_{i=1}^{n}\lambda_{i}(A_{i}+\mu^{-1}\Id)^{-1}]^{-1}\sum_{i=1}^{n}
\lambda_{i}(A_{i}+\mu^{-1}\Id)^{-1}b_{i}}+\\
& q_{[\sum_{i=1}^{n}\lambda_{i}(A_{i}+\mu^{-1}\Id)^{-1}]^{-1}}(\sum_{i=1}^{n}
\lambda_{i}(A_{i}+\mu^{-1}\Id)^{-1}b_{i}),
\end{align*}
we obtain that $\pflm(x^*)=$
\begin{align*}
 &
q_{[{\sum_{i=1}^{n}\lambda_{i}(A_{i}+
\mu^{-1}\Id)^{-1}}]^{-1}-\mu^{-1}\Id}(x^*)
+\scal{x^*}{[\sum_{i=1}^{n}\lambda_{i}(A_{i}+\mu^{-1}\Id)^{-1}]^{-1}\sum_{i=1}^{n}
\lambda_{i}(A_{i}+\mu^{-1}\Id)^{-1}b_{i}}+
\\
& q_{[\sum_{i=1}^{n}\lambda_{i}(A_{i}+\mu^{-1}\Id)^{-1}]^{-1}}(\sum_{i=1}^{n}
\lambda_{i}(A_{i}+\mu^{-1}\Id)^{-1}b_{i})-
\sum_{i=1}^{n}\lambda_{i}\big(q_{(A_{i}+\mu^{-1}\Id)^{-1}}(b_{i})-
r_{i}\big),
\end{align*}
which is \eqref{quadratic}. The remaining claims are immediate from \eqref{quadratic} and
that $\average=0$ when $(\forall\ i)\ A_{i}=0$ by Proposition~\ref{sameoperator}.
\end{proof}

We are ready for our main result:

\begin{theorem}[harmonic-resolvent-arithmetic average inequality and
limits]\label{matrixequ} ~\\
Let $A_{1},\ldots,A_{n}\in \SNPP$. We have
\begin{enumerate}
\item
\begin{equation}\label{harmarithmetic}
\harm\preceq \average \preceq \arithmetic;
\end{equation}
In particular, $\average \in \SNPP$.
\item $\average\rightarrow \arithmetic$ when $\mu\rightarrow 0^{+}$.
\item $\average\rightarrow\harm$ when $\mu\rightarrow +\infty$.
\end{enumerate}
\end{theorem}
\begin{proof}
(i). According to Fact~\ref{t:ineqs},
\begin{equation}
\label{key}
(\lambda_{1}f_{1}^*+\cdots +\lambda_{n}f_{n}^*)^{*}\leq\pflm\leq \lambda_{1}f_{1}+\cdots+\lambda_{n}f_{n}.
\end{equation}
Let $f_{i}=q_{A_{i}}$. Using $(q_{A_{i}})^*=q_{A_{i}^{-1}}$ (by Fact~\ref{qconjugate})
and Lemma~\ref{quadraticprox} we have
\begin{align}
(\lambda_{1}f_{1}^*+\cdots +\lambda_{n}f_{n}^*)^{*}&=(\lambda_{1}q_{A^{-1}_{1}}+\cdots+\lambda_{n}q_{A^{-1}_{n}})^*
=(q_{\lambda_{1}A_{1}^{-1}+\cdots+\lambda_{n}A_{n}^{-1}} )^{*}\nonumber\\
& =q_{(\lambda_{1}A_{1}^{-1}+\cdots+\lambda_{n}A_{n}^{-1})^{-1}}
=q_{\harm}.\label{thanksgiving1}\\
\lambda_{1}f_{1}+\cdots+\lambda_{n}f_{n} & =q_{\lambda_{1}A_{1}+\cdots +\lambda_{n}A_{n}}=q_{\arithmetic},\label{thanksgiving2}\\
\pflm &=q_{\average}. \label{thanksgiving3}
\end{align}
Then \eqref{key} becomes
$$q_{\harm}\leq q_{\average}\leq q_{\arithmetic}.$$
As $q_{X}\leq q_{Y} \Leftrightarrow X\preceq Y$, \eqref{harmarithmetic} is established.
Since $A_{i}\in \SNPP, A_{i}^{-1}\in \SNPP, \lambda_{1}A_{1}^{-1}+\cdots
+\lambda_{n}A_{n}^{-1}\in \SNPP$,
we have $\harm=(\lambda_{1}A_{1}^{-1}+\cdots
+\lambda_{n}A_{n}^{-1})^{-1}\in \SNPP$, thus
$\average\in \SNPP$ by \eqref{harmarithmetic}. (Alternatively,
apply Theorem~\ref{t:blabla}.)

(ii) and (iii):
Observe that ($\forall i$)
$\big(\lambda_i\timess f_i)^*=\lambda_{i}f_{i}^*=\lambda_{i}q_{A_{i}^{-1}}$ has full domain, by \cite[Theorem~16.4]{Rocky},
$$(\lambda_{1}f_{1}^*+\cdots +\lambda_{n}f_{n}^*)^{*}=\big(\lambda_1\timess f_1\Box
\cdots\Box\lambda_n\timess f_n\big).$$
By Fact~\ref{t:mono}, $\forall x\in\RR^N$ one has
$$\lim_{\mu\rightarrow 0^{+}}\pflm(x)=(\lambda_{1}f_{1}+\cdots +\lambda_{n}f_{n})(x),$$
$$\lim_{\mu\rightarrow +\infty}\pflm(x)=(\lambda_{1}f_{1}^*+\cdots+\lambda_{n}f_{n}^{*})^*(x).$$
Since $(\forall i) \ f_{i}, f_{i}^*$
are differentiable on $\RR^N$,
so is $\pflm$ by \cite[Corollary~7.7]{bglw}.
According to Fact~\ref{rockstar}, $\forall x$
\begin{equation}\label{martial1}
\lim_{\mu\rightarrow 0^{+}}\nabla \pflm(x)=\lambda_{1}\nabla f_{1}(x)+\cdots +\lambda_{n} \nabla f_{n}(x),
\end{equation}
\begin{equation}\label{martial2}
\lim_{\mu\rightarrow +\infty}\nabla\pflm(x)=\nabla (\lambda_{1}f_{1}^*+\cdots+\lambda_{n}f_{n}^{*})^*(x).
\end{equation}
Moreover, the convergences in \eqref{martial1}-\eqref{martial2} are uniform on every closed bounded subset of
$\RR^N$.
Now it follows from ~\eqref{thanksgiving1}-\eqref{thanksgiving3} that
$\nabla\pflm=\average$, $\nabla (\lambda_{1} f_{1}+\cdots +\lambda_{n} f_{n})=\arithmetic$,
$\nabla (\lambda_{1}f_{1}^*+\cdots+\lambda_{n}f_{n}^{*})^*
=\harm.$
\eqref{martial1}-\eqref{martial2} transpire to
\begin{equation}\label{julian1}
\lim_{\mu\rightarrow 0^{+}}\average x =\arithmetic x,
\end{equation}
\begin{equation}\label{julian2}
\lim_{\mu\rightarrow + \infty}\average x=\harm x,
\end{equation}
where the convergences are uniform on every closed bounded subset of
$\RR^N$.
Hence (ii) and (iii) follow from \eqref{julian1} and \eqref{julian2}.
\end{proof}

Note that in Theorem~\ref{matrixequ}(ii),(iii), there is no ambiguity since all norms in finite dimensional spaces
are equivalent.

\begin{definition} A function $g\colon \mathbb{D}\to\SN$,
where $\mathbb{D}$ is a convex subset of $\SN$, is matrix convex if
$\forall A_{1}, A_{2}\in \mathbb{D}$, $\forall \lambda\in [0,1]$,
$$g(\lambda A_1+(1-\lambda) A_{2}) \preceq\lambda g(A_{1})+(1-\lambda)g(A_{2}).$$
Matrix concave functions are defined similarly.
\end{definition}
It is easy to see that a symmetric matrix valued function $g$ is matrix concave (resp. convex) if and only if
$\forall x\in \RR^N$ the function
$A\mapsto q_{g(A)}(x)$ is concave (resp. convex).
Some immediate consequences of Theorem~\ref{matrixequ} on matrix-valued functions are:
\begin{corollary}\label{known1}
Assume that $(\forall i)\ A_{i}\in \SNPP$ and $\sum_{i=1}^{n}\lambda_{i}=1$ with $\lambda_{i}>0$. Then
$$(\lambda_{1}A_{1}+\cdots +\lambda_{n}A_{n})^{-1}\preceq \lambda_{1}
 A_{1}^{-1}+\cdots+\lambda_{n}A_{n}^{-1}.$$
Consequently, the matrix function $X\mapsto X^{-1}$
is matrix convex on $\SNPP$.
\end{corollary}
\begin{proof}
Apply Theorem~\ref{matrixequ} equation~\eqref{harmarithmetic} for
$\bA =(A_{1}^{-1},\cdots, A_{n}^{-1}).$
\end{proof}

\begin{corollary}\label{known2}
For every $\mu>0$, the resolvent average matrix function
$\bA\mapsto \average$ given by
\begin{align}\label{amanda}
&(A_{1},\cdots, A_{n})\mapsto [\lambda_{1}(A_{1}+\mu^{-1}\Id)^{-1}+\cdots+\lambda_{n}(A_{n}+\mu^{-1}\Id)^{-1}]^{-1}-\mu^{-1}\Id \nonumber \\
& \mbox{ is matrix concave on $\SNPP\times\cdots\times \SNPP$.}
\end{align}

For each $\lambda=(\lambda_{1}, \cdots, \lambda_{n})$ with $\sum_{i=1}^{n}\lambda_{i}=1$ and
$\lambda_{i}>0\ \forall i$,
the harmonic average matrix function
\begin{equation}\label{known1.5}
(A_{1},\cdots, A_{n})\mapsto (\lambda_{1}A_{1}^{-1}+\cdots+\lambda_{n}A_{n}^{-1})^{-1} \mbox{ is matrix concave}
\end{equation}
 on $\SNPP\times\cdots\times \SNPP$.
Consequently, the harmonic average function
\begin{equation}\label{known3}
(x_{1},\cdots, x_{n})\mapsto \frac{1}{x_{1}^{-1}+\cdots +x_{n}^{-1}} \mbox{ is concave}
\end{equation}
 on $\RR_{++}\times\cdots\times\RR_{++}$.
\end{corollary}
\begin{proof} Set $f_{i}=q_{A_{i}}$.
Then $\forall x\in \RR^N$, we have $\pflm(x)=$
$$\min_{\lambda_{1}x_{1}+\cdots +\lambda_{n}x_{n}=x}\bigg(\big(\lambda_{1}q_{A_{1}}(x_{1})+
\cdots +\lambda_{n}q_{A_{n}}(x_{n})\big)+\big(\mu^{-1}\lambda_{1}q_{\Id}(x_{1})+\cdots+\mu^{-1}\lambda_{n}
q_{\Id}(x_{n})\big)\bigg)
-\mu^{-1}q_{\Id}(x).$$
Since for each fixed $(x_{1},\ldots,x_{n})$,
$$(A_{1},\cdots, A_{n})\mapsto
\big(\lambda_{1}q_{A_{1}}(x_{1})+
\cdots +\lambda_{n}q_{A_{n}}(x_{n})\big)+\big(\mu^{-1}q_{\Id}(x_{1})+\cdots+\mu^{-1}q_{\Id}(x_{n})\big)
,$$
is affine, being the infimum of affine functions we have that $\forall x$ the function
$$(A_{1},\cdots, A_{n})\mapsto \pflm(x),$$
is concave. As
$\pflm(x)=q_{\average}(x)$ by Lemma~\ref{quadraticprox},
this shows that $\forall x\in\RR^N$ the function
$$\bA=(A_{1},\cdots,A_{n})\mapsto q_{\average}(x) \mbox{ is concave,}$$ so
$\bA \mapsto \average$
is matrix concave.

Now by Theorem~\ref{matrixequ}(iii), $\average\rightarrow\harm$ when $\mu\rightarrow +\infty$.
This and \eqref{amanda}
implies that
$$\bA\mapsto\harm,$$
is also matrix concave. \eqref{known3} follows from \eqref{known1.5}
by setting $N=1$ and $\lambda_{1}=\cdots=\lambda_{n}=1/n$.
\end{proof}
\begin{remark} Corollary~\ref{known1} is well-known,
cf.\ \cite[Proposition~2.56 on page~73]{RockWets}.
Corollary~\ref{known2} \eqref{known3} is also well-known,
cf.\ \cite[Exercise~3.17 on page~116]{boyd}.
\end{remark}

We proceed to show that resolvent averages of matrices enjoy self-duality.
\begin{theorem}[self-duality]\label{mainresult}
 Let $(\forall \ i)\ A_{i}\in \SNPP$ and $\mu>0$. Assume that
$\sum_{i=1}^{n}\lambda_{i}=1$ with $\lambda_{i}>0$. Then
\begin{equation}\label{niceselfdual}
\big[\average\big]^{-1}=\averageinverse, \mbox{ i.e.,}
\end{equation}
\begin{align*}
&\bigg[\bigg(\lambda_{1}(A_1+\mu^{-1}\Id)^{-1}+\cdots +\lambda_{n}(A_{n}+\mu^{-1}\Id)^{-1}\bigg)^{-1}-\mu^{-1}\Id\bigg]^{-1}=\\
& \bigg(\lambda_{1}(A_1^{-1}+\mu\Id)^{-1}+\cdots +\lambda_{n}(A_{n}^{-1}+\mu\Id)^{-1}\bigg)^{-1}-\mu\Id.
\end{align*}
In particular, for $\mu=1$, $[\averageonelambda]^{-1}=\averageoneinverse.$
\end{theorem}

\begin{proof} Let $f_{i}=q_{A_{i}}$.
By Fact~\ref{t:Fenchel},
$\big(\pflm\big)^* = p_{\mu^{-1}}({\fettf^*,\fettla}),$
taking subgradients both sides, followed by using Fact~\ref{conjgrad}, we obtain that
$$\partial \big(\pflm\big)^*= \big(\partial \pflm\big)^{-1}=\partial \big(p_{\mu^{-1}}({\fettf^*,\fettla})\big).$$
By Lemma~\ref{quadraticprox}, $\pflm=q_{\average}, p_{\mu^{-1}}({\fettf^*,\fettla})=q_{\averageinverse},$ we have
$$\partial \pflm = \average, $$
$$\partial p_{\mu^{-1}}({\fettf^*,\fettla})=\averageinverse.$$
Hence
$$\big[\average\big]^{-1}=\averageinverse,$$
as claimed.
\end{proof}
\begin{remark} Although the harmonic and arithmetic average
lack self-duality, they are dual to each other:
$$[\harm]^{-1}=\lambda_{1}A_{1}^{-1}+\cdots+\lambda_{n}A_{n}^{-1}=\arithmeticinverse,$$
$$[\arithmetic]^{-1}=\big[\lambda_{1}(A_{1}^{-1})^{-1}+\cdots+\lambda_{n}(A_{n}^{-1})^{-1}\big]^{-1}
=\harminverse.$$
\end{remark}

\section{A comparison to weighted geometric means}\label{compare}
If $A,B\in \SNPP$, the geometric mean is defined by
$$A\sharp B=A^{1/2}\big(A^{-1/2}BA^{-1/2}\big)^{1/2}A^{1/2}.$$
In general, the geometric mean of $A_{1},\ldots, A_{n}\in \SNP$ for $n\geq 3$ is defined either as the limit of an inductive procedure or by the Riemannian distance without a closed form \cite{ando,petz,moakher,kubo}.
To compare the resolvent average with the well-known geometric mean, we restrict our attention to
non-negative real numbers ($1\times 1$ matrices).
When $\bA=\fettx =(x_{1},\ldots, x_{n})$ with $x_{i}\in\RR_{+}$ and $\mu=1$, we write
$$\averagex=\average=\big(\lambda_{1}(x_{1}+1)^{-1}+\cdots+\lambda_{n}(x_{n}+1)^{-1}\big)^{-1}-1,
$$
and $\fettx^{-1}=(1/x_{1},\ldots, 1/x_{n})$ when $(\forall i)\ x_{i}\in\RR_{++}$.
\begin{proposition} Let $(\forall\ i)\ x_{i}>0, y_{i}>0$. We have
\begin{enumerate}
\item[(i)]\emph{({\bf harmonic-resolvent-arithmetic mean inequality})}:
\begin{equation}\label{hram}
\big(\lambda_{1}x_{1}^{-1}+\cdots+\lambda_{n} x_{n}^{-1}\big)^{-1}\leq
\averagex\leq \lambda_{1}x_{1}+\cdots +\lambda_{n}x_{n}.
\end{equation}
Moreover, $\averagex=\lambda_{1}x_{1}+\cdots+\lambda_{n}x_{n}$ if and only if $x_{1}=\cdots=x_{n}$.
\item[(ii)]\emph{({\bf self-duality})}: $[\averagex]^{-1}=\averageinversex$.
\item[(iii)] If $\fettx=(x_{1},\ldots,x_{1})$, then $\averagex=x_{1}$.
\item[(iv)]If $\fettx=(x_{1},x_{1}^{-1}, x_{2},x_{2}^{-1},\ldots,
x_{n},x_{n}^{-1})$ and $\fettla=(\tfrac{1}{2n},\ldots,\tfrac{1}{2n})$,
then $\averagex=1$.
\item[(v)] The function $\fettx\mapsto \averagex$ is concave on $\RR_{++}\times\cdots\times\RR_{++}$.
\item[(vi)] If $\fettx\succeq \fetty$, then $\averagex\geq \averagey$.
\end{enumerate}
\end{proposition}
\begin{proof}
(i): For \eqref{hram}, apply Theorem~\ref{matrixequ}(i) with $\mu=1$. Now
 $\averagex=\lambda_{1}x_{1}+\cdots+\lambda_{n}x_{n}$ is equivalent to
 \begin{equation}\label{newform}
 \big(\lambda_{1}(x_{1}+1)^{-1}+\cdots+\lambda_{n}(x_{n}+1)^{-1}\big)^{-1}=\lambda_{1}x_{1}+\cdots+\lambda_{n}x_{n}+1,
 \end{equation}
As $\sum_{i=1}^{n}\lambda_{i}=1$, \eqref{newform} is the same as
 $$\lambda_{1}\frac{1}{(x_{1}+1)}+\cdots+\lambda_{n}\frac{1}{(x_{n}+1)}=\frac{1}{\lambda_{1}(x_{1}+1)+\cdots+
 \lambda_{n}(x_{n}+1)}.$$
 Since the function $x\mapsto 1/x$ is strictly convex on
 $\RR_{++}$, we must have $x_{1}=\cdots=x_{n}$.

 (ii): Theorem~\ref{mainresult}. (iii): Proposition~\ref{sameoperator}.
  (iv): Proposition~\ref{primalinverse}. (v): Corollary~\ref{known2}. (vi): Proposition~\ref{monotonecompare}.
\end{proof}

 Recall the \emph{weighted geometric mean}:
$$\geoaverage=x_{1}^{\lambda_{1}}x_{2}^{\lambda_{2}}\cdots x_{n}^{\lambda_{n}}.$$
$\geoaverage$ always has the following properties:
\begin{fact} Let $(\forall\ i)\ x_{i}>0, y_{i}>0$. We have
\begin{enumerate}
\item[(i)] \emph{({\bf harmonic-geometric-arithmetic mean inequality})}:
$$\big(\lambda_{1}x_{1}^{-1}+\cdots+\lambda_{n} x_{n}^{-1}\big)^{-1}\leq
\geoaverage\leq \lambda_{1}x_{1}+\cdots +\lambda_{n}x_{n}.$$
Moreover, $\geoaverage=\lambda_{1}x_{1}+\cdots +\lambda_{n}x_{n}$ if and only $x_{1}=\cdots=x_{n}$.
\item[(ii)]\emph{({\bf self-duality})}: $[\geoaverage]^{-1}=\geoinverse$.
\item[(iii)] If $\fettx=(x_{1},\ldots,x_{1})$, then $\geoaverage=x_{1}$.
\item[(iv)] If $\fettx=(x_{1},x_{1}^{-1}, x_{2},x_{2}^{-1},\ldots,
x_{n},x_{n}^{-1})$ and $\fettla=(\tfrac{1}{2n},\ldots, \tfrac{1}{2n})$, then $\geoaverage=1$.
\item[(v)] The function $\fettx\mapsto \geoaverage$ is concave on $\RR_{++}\times\cdots\times\RR_{++}$.
\item[(vi)] If $\fettx\succeq \fetty$, then $\geoaverage\geq \geoaveragey$.
\end{enumerate}
\end{fact}
\begin{proof} (i): See \cite[page 29]{Rocky}. (ii)-(iv) and (vi) are simple. (v): See \cite[Example 2.53]{RockWets}.
\end{proof}

The means $\averagex$ and $\geoaverage$ have strikingly
similar properties. Are they the same?
\begin{example}
(i). Let $\lambda=(\thalb,\thalb)$. When $x=(0,1)$, $\geoaverage=0$ but
$\averagex=\tfrac{1}{3}$, so $\geoaverage
\neq\averagex$.

(ii). \emph{Is it right that $\geoaverage\leq \averagex$
$\forall\ x\in \RR_{++}^2$?}
The answer is also no.
Assume to the contrary that
$\geoaverage\leq \averagex, \ \forall \ \fettx \in \RR_{++}\times \RR_{++}$. Taking inverse both sides, followed by applying the self-duality
of $\geoaverage, \averagex$, gives
$$\geoaverage^{-1}\geq \averagex^{-1}=\averageinversex\geq \geoinverse=\geoaverage^{-1},$$
and this gives that $\geoaverage^{-1}=\averagex^{-1}$ so that $\geoaverage=\averagex$.
This is a contradiction to (i).
\end{example}

Finally, we note that the resolvent average can be defined for general monotone operators and that
Theorem~\ref{mainresult} holds even when $A_1,\ldots,A_n$
are monotone operators
(not necessarily positive semi-definite matrices), in that situation one needs to use \emph{set-valued
inverses}.
This and further details on the resolvent average for
general monotone operators will appear in the forthcoming paper
\cite{bws}.

\section*{Acknowledgment}
Heinz Bauschke was partially supported by the Canada Research
Chair program and by the Natural Sciences and Engineering Research
Council of Canada.
Sarah Moffat was partially supported by the Natural Sciences and
Engineering Research Council of Canada.
Xianfu Wang was partially supported by the Natural Sciences and
Engineering Research Council of Canada.

\end{document}